\documentclass[a4paper,12pt]{article}
\usepackage{amsfonts}
\usepackage{}
\usepackage{mathrsfs}
\usepackage{algorithm, algpseudocode}
\usepackage{amsmath,amssymb,amsthm,latexsym,amsfonts}
\usepackage{pstricks}
\usepackage{enumerate}

\usepackage{graphicx}
\usepackage{color}
\usepackage{ifpdf}

\usepackage{caption}

\begin{document}\baselineskip 0.8cm

\setlength{\oddsidemargin}{0in} \setlength{\textwidth}{6.5in}
\setlength{\textheight}{8.5in} \setlength{\topmargin}{-.5in}
\def\BBox{\vrule height 0.5em width 0.6em depth 0em}

\title{Constructions of graphs and trees with partially prescribed spectrum\footnote{Supported by NSFC No.11371205 and 11531011, and
PCSIRT.}}
\author{Xueliang Li\thanks{Center for Combinatorics and LPMC, Nankai University, Tianjin 300071, China.
email:lxl@nankai.edu.cn}
\and
Wasin So\thanks{Department of Mathematics, San Jose State University, San Jose, California
95192-0103, USA. email: wasin.so@sjsu.edu}
\and
Ivan Gutman\thanks{Faculty of Science, University of Kragujevac, P.O. Box 60, 34000 Kragujevac, Serbia.
email:gutman@kg.ac.rs}
}
\maketitle

It is shown how a connected graph and a tree with partially prescribed spectrum
can be constructed. These constructions are based on a recent result of
Salez that every totally real algebraic integer is an eigenvalue of a tree.
Our result implies that for any (not necessarily connected) graph $G$, there is a tree $T$ such that the characteristic
polynomial $P(G,x)$ of $G$ can divide the characteristic polynomial $P(T,x)$ of $T$, i.e., $P(G,x)$ is
a divisor of $P(T,x)$.\\[2mm]

\section{Introduction}

Graph eigenvalues have been studied intensively \cite{brouwer_haemers, cvetk, cvetkovic_doob_sachs}, and
they are very special real numbers. Indeed,  they are roots of monic
integral polynomials with only real roots, i.e., they are totally
real algebraic integers. It is natural for one to wonder whether the
converse is true. Forty years ago, Hoffman \cite{hoffman} conjectured that this is
true, which eventually was confirmed by Estes in 1992 \cite{estes}.
\begin{description}
\item[Theorem 1.1.] {\upshape \cite{estes}}
Every totally real algebraic integer is an eigenvalue of a (connected) graph.
\end{description}
Recently, Salez \cite{salez} strengthened the result with a simpler proof.
\begin{description}
\item[Theorem 1.2.] {\upshape \cite{salez}}
Every totally real algebraic integer is an eigenvalue of a tree.
\end{description}

The next natural question is which collection of totally real algebraic integers
forms the spectrum of a graph. Of course, there are many more necessary conditions
on such collections. Below, we list just a few.

\begin{description}
\item[Lemma 1.3.] If $S=\{\lambda_1 \ge \lambda_2 \ge \cdots \ge \lambda_n\}$ is the spectrum
of a graph of order $n$, then

\begin{enumerate}
	\item $S$ contains all the conjugates of each $\lambda_i$,
	\item $\lambda_1 + \cdots + \lambda_n =0$,
	\item $\lambda_1^2+\cdots+\lambda_n^2 \le n(n-1)$,
	\item  $\lambda_1 \le n-1$,
	\item $|\lambda_n| \le \lambda_1$. \BBox
\end{enumerate}
\end{description}

Unfortunately, these conditions are far from being sufficient, as the next example shows.

\begin{description}
\item[Example 1.4.] $\{2, 1,-1, -2 \}$ satisfies  all conditions listed in Lemma 1.3,
but it is not the spectrum of any graph of order 4.

\textit{Proof:} Suppose that there is a graph $G$ such that
 $Spec(G)=\{2, 1,-1, -2 \}$. Then $G$ is bipartite because $Spec(G)$ is symmetric about 0.
	Hence the number of edges of $G$ is at most 4, because $G$ is a bipartite graph of order 4.
	On the other hand,  the number of edges of $G$, computed by means of the eigenvalues,
would be   $\frac {1}{2} \big[2^2+1^2 +(-1)^2+(-2)^2 \big]=5$, contradiction! \BBox
\end{description}

The problem of finding necessary and sufficient conditions for a set of totally real algebraic
integers to be the spectrum of a graph seems intractable! Instead, we tackle the following modified
problem.

\begin{description}
\item[Problem 1.5.] Construct a connected graph such that its spectrum contains a given set of
totally real algebraic integers.
\end{description}

In Section 2, we accomplish such a construction via Knonecker product of matrices.
In Section 3, we strengthen the result by constructing a tree via an appropriate
graph operation.

\section{Construction of connected graphs}

Recall some facts about Kronecker product of matrices:

\begin{itemize}
	\item[Fact 1.] $Spec(A \otimes B) = \{ \alpha \beta: \alpha \in Spec(A), \beta \in  Spec(B)\}$
	\item[Fact 2.]  $Spec(A \otimes I + I \otimes B)= \{ \alpha + \beta: \alpha \in Spec(A), \beta \in  Spec(B)\}$
	\item[Fact 3.] If $A$ and $B$ are adjacency matrices, then $A \otimes B$ is also an adjacency matrix.
	\item[Fact 4.] If $A$ and $B$ are adjacency matrices, then $A \otimes I + I \otimes B$ is also an adjacency matrix.
\end{itemize}

In view of Facts 3 and 4, we introduce two graph products as follows:

\begin{description}
\item[Definition 2.1.] Given two graphs $G$ and $H$, define a new graph $G+H$ such that
its adjacency matrix is given by $A(G+H)=A(G) \otimes I +I \otimes A(H)$.
\item[Definition 2.2.] Given two graphs $G$ and $H$, define a new graph $G \times H$ such that
its adjacency matrix is given by $A(G \times H)=A(G) \otimes  A(H)$.
\end{description}

Moreover, if $G$ and $H$ are connected, then $G+H$ is also connected. It is well-known that $G \times H$ is connected if and only if both $G$ and $H$ are connected and one of $G$ and $H$ contains a
cycle of odd length, i.e., one of them is non-bipartite.

Using Facts 1 and 2, we have
\[ Spec(G+H) = Spec(G) + Spec(H), \] and
\[ Spec(G \times H) = Spec(G) \cdot Spec(H).\]

Since the path $P_5$ of order 5 has $Spec(P_5)=\{-\sqrt{3}, -1, 0, 1, \sqrt{3}\}$ and
the cycle $C_3$ of order 3 has $Spec(C_3)=\{-1, -1, 2\}$,
we have that the graph $F=P_5 + C_3$ has $0=1+(-1)$ and $1=(-1)+2$ as its eigenvalues. Obviously, $F$ is connected and non-bipartite since it contains an odd cycle $C_3$. As will be seen in the following, we only need non-bipartite graphs $F$ that have 0 and 1 as its eigenvalues. The above says the existence of such graphs. Actually, there are such graphs of small order and size. For example, the graph obtained by attaching two pendant vertices and a 2-vertex path to the same vertex of the triangle. This graph has 7 vertices.

\begin{description}
 \item[Lemma 2.3.]  Given a connected graph $G$ such that $\alpha \in Spec(G)$.
 Then there is a connected graph $H$ such that $ 0, \alpha \in Spec(H)$.

 \textit{Proof:}
 Take $H=F \times G$. Then $H$ is connected since $F$ is non-bipartite. Moreover,
 since $F$ contains 0 and 1 as its eigenvalues, we have $0, \alpha \in Spec(H)$  \BBox

 \item[Theorem 2.4.] Let $\alpha_1, \ldots, \alpha_p$ be totally real algebraic integers. Then
 there is a connected graph $H$ such that $\{\alpha_1, \ldots,\alpha_p\} \subseteq Spec(H)$.

 \textit{Proof:} We prove, by induction on $p$, a stronger statement:
 there is a connected graph $H$ such that $\{0, \alpha_1, \ldots,\alpha_p\} \subseteq Spec(H)$.

 Consider $p=1$. By Theorem 1.1, there is a graph $G$ such that $\alpha_1 \in Spec(G)$.
 Without loss of generality, we can assume that $G$ is connected. Now, by Lemma 2.3,
 there is a connected graph $H$ such that $ 0, \alpha_1 \in Spec(H)$.

 Consider $p>1$. By the induction assumption, there is a connected graph $K$ such that
 $\{0, \alpha_1, \ldots,\alpha_{p-1}\} \subseteq Spec(K)$. Applying the case $p=1$,
 we have a connected graph $G$ such that $0, \alpha_p \in Spec(G)$. Take $H=K+G$.
 Then $H$ is connected because both $K$ and $G$ are connected. Moreover,
\[ 0, \alpha_1, \ldots, \alpha_{p-1}, \alpha_p \in \{0, \alpha_1, \ldots,\alpha_{p-1}\}+\{0, \alpha_p \}
 \subseteq Spec(K)+ Spec(G) = Spec(H). \BBox \]

 \end{description}

\section{Construction of trees}

\begin{description}
\item[Lemma 3.1:] Let $A$ and $B$ be square matrices. Then
\[ Spec \left( \left[ \begin{array}{ccc} A & F &F\\E&B&0\\E&0&B \end{array} \right] \right) = Spec( B ) \bigcup
Spec \left( \left[ \begin{array}{cc} A & 2F \\E&B  \end{array} \right] \right).\]

\textit{Proof:} Note that
\[ \left[ \begin{array}{ccc} I & 0 &0\\0&I&0\\0&I&I \end{array} \right] =
\left[ \begin{array}{ccc} I & 0 &0\\0&I&0\\0&-I&I \end{array} \right]^{-1}.\]
Then the following matrix identity is in fact a similarity transformation:
\[ \left[ \begin{array}{ccc} I & 0 &0\\0&I&0\\0&-I&I \end{array} \right]
\left[ \begin{array}{ccc} A & F &F\\E&B&0\\E&0&B \end{array} \right]
\left[ \begin{array}{ccc} I & 0 &0\\0&I&0\\0&I&I \end{array} \right]
= \left[ \begin{array}{ccc} A & 2F &F\\E&B&0\\0&0&B \end{array} \right].\]
Hence $\left[ \begin{array}{ccc} A & F &F\\E&B&0\\E&0&B \end{array} \right]$
 and $\left[ \begin{array}{ccc} A & 2F &F\\E&B&0\\0&0&B \end{array} \right]$
 have the same spectrum, and so the conclusion follows. \BBox
\end{description}

Given disjoint graphs $G$, $H_i$, and $H_i'$  such that
$H_i$ and $H_i'$ are isomorphic for $i=1,2,\ldots,p$.
Let $x_i \,,\,i=1,2,\ldots,p$, be vertices of $G$ (not necessarily different).
Let $v_i$ be a vertex of $H_i$, and $v_i'$ a vertex of $H_i'$.
Construct a graph $G \circ [H_1, \cdots, H_p]$ by connecting
$x_i$ to both $v_i$ and $v_i'$ with new edges, for  $i=1,2,\ldots,p$.

\begin{description}
\item[Lemma 3.2.] $Spec(H_1 \cup \cdots \cup H_p) \subseteq Spec(G \circ [H_1, \cdots, H_p])$.

\textit{Proof:} Let $H=H_1 \cup \cdots \cup H_p$
and  $H'=H_1' \cup \cdots \cup H_p'$. Since $H_i$ and $H_i'$ are isomorphic,
 $H$ and $H'$ are also isomorphic.
Hence, by a suitable labeling, we have $A(H)=A(H')$ and
\[ A(G \circ [H_1, \cdots, H_p]) = \left[ \begin{array}{ccc} A(G) & E^T &E^T\\E&A(H)&0\\E&0&A(H) \end{array} \right]. \]
Consequently, by Lemma 3.1,
\begin{eqnarray*}
 Spec(H_1 \cup \cdots \cup H_p) & = & Spec(A(H)) \\
 & \subseteq&  Spec(A(G \circ [H_1, \cdots, H_p])) \\
 &=& Spec(G \circ [H_1, \cdots, H_p])\,. \;\;\;  \BBox
 \end{eqnarray*}

 \item[Theorem 3.3.] Let $\alpha_1, \ldots, \alpha_p$ be totally real algebraic integers. Then
 there is a tree $T$ such that $\{\alpha_1, \ldots,\alpha_p\} \subseteq Spec(T)$.

 \textit{Proof:}
 For each   totally real algebraic integer $\alpha_i$, by Theorem 1.2, there
 is a tree $T_i$ whose spectrum contains $\alpha_i$. Take $G$ to be any tree (say, just a singleton).
 By Lemma 3.2,
 $Spec(T_1 \cup \cdots \cup T_p) \subseteq Spec(G \circ [T_1, \cdots, T_p])$
 and so \[ \{\alpha_1, \ldots,\alpha_p\} \subseteq Spec(G \circ [T_1, \cdots, T_p]). \]
 Moreover, $T=G \circ [T_1, \cdots, T_p]$ is a tree because $G$ and $T_i$ are all trees. \BBox

 \item[Example 3.4.] Note that $Spec(K_2)=\{-1, 1\}$, and $Spec(K_{1,4})=\{-2,0, 0, 0, 2\}$.
 Hence, by the construction in the proof of Theorem 3.4, $E_1 \circ[K_2, K_{1,4}]$ is a tree
  whose spectrum contains $\{-2, -1, 1, 2\}$.

Let $G$ be a graph. A $k$-matching of $G$ is a set of $k$ edges such that any two distinct edges in the set do not share a common vertex.
The matching polynomial $m(G,x)$ of $G$ is defined as
$$
m(G,x)=\sum_{k\geq 0}(-1)^k m(G,k)x^{n-2k},
$$
where $m(G,k)$ denotes the number of $k$-matchings in $G$, and $m(G,0)=1$ by convention.

For matching polynomials, we know from \cite{cvetk} that
for any (not necessarily connected) graph $G$, all the roots of $m(G,x)$ are totally real algebraic integers, and moreover, there is
a tree $T$ such that $m(G,x)$ is a divisor of $m(T,x)$. The next result says that the same thing holds for characteristic polynomials of graphs.

 \item[Corollary 3.5.] For any (not necessarily connected) graph $G$, there is a tree $T$ such that the characteristic
 polynomial $P(G,x)$ of $G$ can divide the characteristic polynomial $P(T,x)$ of $T$, i.e., $P(G,x)$ is
 a divisor of $P(T,x)$.

 \textit{Proof:}
 Since all the roots of $P(G,x)$ are totally real algebraic integers, by Theorem 3.3 there
 is a tree $T$ whose spectrum contains all the roots of $P(G,x)$, and hence the conclusion
 follows. \BBox

 A real polynomial is unimodal if the sequence of the coefficients of the polynomial is unimodal, 
 i.e., first increasing, and then decreasing, with only one peak.

 \item[Corollary 3.6.] For any totally real algebraic polynomial $f(x)$, there is another
 totally real algebraic polynomial $g(x)$ such that $f(x)g(x)$ is unimodal.

 \textit{Proof:} From Theorem 3.3, we know that $f(x)$ is a divisor of the characteristic polynomial
 of a tree. It is well-known that the characteristic polynomial of any bipartite graph, and therefore,
 any tree, is unimodal. The conclusion follows immediately.  \BBox

\end{description}

This result means that any totally real algebraic polynomial can be unimodalized. For example, the
characteristic polynomial of an arbitrary graph is usually not unimodal, but it can be unimodalized by another
totally real algebraic polynomial. It could be an interesting question to think about how to unimodalize a
totally real algebraic polynomial by using another totally real algebraic polynomial with a degree as small
as possible.


\begin{thebibliography}{20}

\bibitem{brouwer_haemers}
A.E. Brouwer and W.H. Haemers,
\textit{Spectra of Graphs}, Springer, New York, 2012

\bibitem {cvetk} D.M. Cvetkovi\'c,	M. Doob, I. Gutman,	A. Torgasev,
\textit{Recent Results in the Theory of Graph Spectra}, Ann. Discrete Math. No.36, North-Hollands, 1988.

\bibitem{cvetkovic_doob_sachs}
D.M. Cvetkovi\'c, M. Doob and H. Sachs,
\textit{Spectra of Graphs: Theory and Applications, 3rd edition,}
Johann Ambrosius Barth, Heidelberg, 1995.

\bibitem{estes}
D.R. Estes,
Eigenvalues of symmetric integer matrices,
\textit{J. Number Theory} 42(3) (1992) 292--296.

\bibitem{hoffman}
A.J. Hoffman, Eigenvalues of graphs,
in: D. R. Fulkerson (Ed.), \textit{Studies in Graph Theory, Part II},
Math. Assoc. Amer., Washington, DC, 1975, pp. 225--245.


\bibitem{salez}
J. Salez,
Every totally real algebraic integer is a tree eigenvalue,
\textit{J. Comb. Theory, Series B} 111 (2015) 249--256.


\end{thebibliography}
\end{document}